\documentclass{amsart}
\usepackage{amssymb}
\newtheorem{theorem}{Theorem}[section]

\begin{document}
\begin{center}{\Large Contarctive properties of multifunctions related to uniformity}\\\smallskip
{\normalsize Pratip Chakraborty\footnote{\copyright\,My pleasure if you want to copy}\\
University of Kaiserslautern\\Kaiserslautern, DE
67663\\chakrabo@mathematik.uni-kl.de
 \\ 28/01/2005}

\end{center}\bigskip
\section{Introduction}
In this paper I will try to prove some contactive properties for
multifunctions in uniformly complete metric spaces which are
orbittally complete.
\subsection{Basic definitions and notations}

Let $(X,U)$ be uniform space. A family $\{d_i:i\in I\}$ of
pseudometrics on $X$ with indexing set $I$, Is called an associated
family for the uniformity $U$ if the family
$$\beta =\{V(i,\varepsilon):i\in I;\varepsilon > o\}$$

where
$$ v(i,\varepsilon)=\{(x,y):x,y\in A, i\in I\}$$
is a sub-base for the uniformity$U$. We may assume that $\beta$
itself is a base by adjoining finite intersection members of $\beta$
, if necessary.The corresponding family of psedometrics is called an
associated family for $U$. An associated family for $U$ will be
denoted by $p^*$.Let $A$ be a nonempty subset of a uniform space$X$.
Define
$$D^*(A)=sup\{d_i(x,y):x,y \in A,i\in I\}$$
where
$$\{d_i:i\in I\}=p^*$$
Then $D^*$ is called an augumented diameter of A. Furthermore A is
said to be {\em $p^*$-bounded} if $D^*(A) < \infty$. Let
\[2^X=\{A:
\mbox{A is nonempty, closed and $p^*$-bounded subset of X}\}\] For
any non empty subsets of A and B of X,
define$$d_i(x,A)=inf\{d_i(x,A):a\in A\}, i\in I$$
$$H_i(A,B)=max\{\sup d_i(a,B), \sup d_i(A,b)\}=sup\{|d_i(x,A)-d_i(x,B)|\}$$
where $a\in A,b\in B$ \text{and} $x\in X $
\newline
It is well known that on $2^X $, $H_i$ is a pseudometric, called the
induced Hausdorff pseudometric induced by $d_i,i\in I.$
\newline
Let $(X,U)$ be a unifrom space with an augumented associated family
$p^*.$ $p^*$ also induces a uniformity $U^*$ on $2^X $ defined by
the base $$\beta^*=\{V^*(i,\varepsilon):i\in I,\varepsilon > 0 \}$$
where$$V^*(i,\varepsilon)=\{(A,B):A,B\in 2^X,H_i(A,B)<
\varepsilon\}.$$ The space $(2^X,U^*)$ is a uniform space called the
hyperspace of $(X,U)$.
\newline Let us review two following definitions.
\newline
\large Subfilter :
\newline
\normalsize The collection of alll filters on a given set $x$ is
denoted by $\Phi(X).$ An order relation \newpage is defined on
$\Phi(X)$ by the rule $\mathfrak F_1<\mathfrak F_2$ iff $\mathfrak
F_1 \supset \mathfrak F_2.$ If $\mathfrak F^*<\mathfrak F$, then
$\mathfrak F^*$ is called a subfilter of $\mathfrak F.$
\newline
\large Orbital Completeness : \normalsize Let $(X,U)$ be auniform
space defined by $\{d_i:i\in I\}=p^*.$ if $F:X \rightarrow 2^X$ ia a
multifunction, then \newline \begin{enumerate}
\item  $x \in X$ iscalled a fixed point of $F$ if $\in Fx;$ complete
\item  an orbit of $F$ at a point $x_0 \in X$ is a sequence $\{x_n\}$ given by
$$O(F,x_0)=\{x_n:x_n \in Fx_n-1,n=1,2,...\};$$
\item A uniform space $X$ is called \emph{$F$-orbitally complete}
if every Cauchy filter which is a subfilter of an orbit of F at each
$x \in X$ converges to a point of $X$.
\end{enumerate}

\large Orbital Continuity : \normalsize Let $(X,U)$ be a uniform
space and let $F:X \rightarrow X$ be a function.A single-valued
function $F$ is \emph{orbitally continuous} if $lim(T^{n_i}x)=u$
implies $lim T(T^{n_i}x)=Tu$ for each $x \in X$.
\section{Contraction theorems}
\begin{theorem} Let $(x,U)$ be a \emph{F-orbitally complete} Hausdorff uniform
space defined by $\{d_i:i\in I\}=p^*.$ and $(2^X,U^*)$ a hyperspace
and let $F:X \rightarrow 2^X$ be a continuous function with $Fx$
compact for each $x \in X$. Assume that
$$min\{H_i(Fx,Fy)^r,d_i(x,Fx)d_i(y,Fy)^{r-1},d_i(y,Fy)^r\}$$ $$+ a_{i} min\{d_i(x,Fy),d_i(y,Fx)\}\le [b_i d_i(x,Fx)+c_i d_i(x,y)]d_i(y,fy)^{r-1}$$
for all $i \in I$ and $x,y \in X,$ where $r \ge 1$ is an integer,
$a_i,$ $b_i,$ $c_i$ are real numbers such that $0<b_i+c_i<1$, then
$F$ has a fixed point.
\end{theorem}
\emph{proof} : Let $x_0$ be an arbitrary point in $X$ and consider
the sequence defined by $$x_1 \in Fx_0,x_2 \in Fx_1,......,x_n \in
Fx_{n-1},...$$ Let us assume that $d_i(x_n,Fx_n)> 0$ for each $i \in
I$ and $n=0,1,2,...$(otherwise for some positive integer $n, x_n \in
Fx_n$ as desired.)\newline Let $\mathcal{U} \in U $ be an arbitrary
entourage. Since $\beta$ is a base for $U,$ there exists
$V(i,\varepsilon)\in \beta$ such that $V(i,\varepsilon) \subseteq
\mathcal{U}.$ Now $y\rightarrow d_i(x_0,y)$ is continuous on the
comapact set $Fx_0$ and this implies that there exists $x_1 \in
Fx_0$ such that $d_i(x_0,x_1)=d_i(x_0,Fx_0).$ Similarly, $Fx_1$ is
compact so there exists $x_2 \in Fx_1$ such
 that $d_i(x_1,x_2)=d_i(x_1,Fx_1).$ Counting, we obtain a sequence
$\{x_n\}$ such that $x_n+1 \in Fx_n$ and
$d_i(x_n,x_n+1)=d_i(x_n,Fx_n).$ For $x=x_{n-1},$ and $y=x_n$ by the
condition of the theorem we have
$$min\{H_i(Fx_{n-1},Fx_n)^r,d_i(x_{n-1},Fx_{n-1})d_i(x_n,Fx_n)^{r-1},d_i(x_n,Fx_n)^r\}$$
$$+ a_{i} min\{d_i(x_{n-1},Fx_n),d_i(y,Fx_{n-1})\}\le [b_i d_i(x_{n-1},Fx_{n-1})$$
$$+c_i d_i(x_{n-1},x_n)]d_i(x_n,Fx_n)^{r-1}$$
or since $d_i(x_n,x_{n-1}) = 0,x_n \in Fx_{n-1}.$ Hence we have
$$min \{d_i(x_n,x_{n+1})^r,d_i(x_{n-1},x_n)d_i(x_n,x_n+1)^{r-1}\}$$
$$\le [b_id_i(x_{n-1},x_n) +
c_id_i(x_{n-1},x_n)]d_i(x_n,x_{n+1})^{r-1}$$ and it follows that
$$min \{d_i(x_n,x_{n+1})^r,d_i(x_{n-1},x_n)d_i(x_n,x_n+1)^{r-1}\}$$ $$\le(b_i+c_i)d_i(x_{n-1},x_n)d_i(x_n,x_{n+1})^{r-1}$$
is not possible (as $0$ $< b_i + c_i <1$), we have
$$d_i(x_n,x_{n+1})^r \le(b_i+c_i)d_i(x_{n-1},x_n)d_i(x_n,x_{n+1})^{r-1}$$
or $$d_i(x_n,x_{n+1})^r \le
k_id_i(x_{n-1},x_n)d_i(x_n,x_{n+1})^{r-1},$$ where $k_i = b_i + c_i,
0< k_i < 1.$ \newline proceeding in this mway we get
\begin{eqnarray*}
d_i(x_n,x_{n+1})&\le&k_id_i(x_{n-1},x_n)\\&\le&k_i^2d_i(x_{n-2},x_{n-1})
\end{eqnarray*}
$$.$$ $$.$$ $$.$$   $$
\le k_i^nd_i(x_0,x_1)$$ Hence we obtain
\begin{eqnarray*}
d_i(x_n,x_m)&\le& d_i(x_n,x_{n+1}) + d_i(x_{n+1},x_{n+2} + ....+
d_i(x_{m-1},x_m)\\&\le& (k_i^n + k_i^{n+1} + ....
+k_i^{m-1})d_i(x_0,x_1)\\&\le&k_i^n(1 + k_i + .... +
k_i^{m-n-1}\\&\le& \frac{k_i^n}{1-k} (1 + k_i + .... +
k_i^{m-n-1})d_i(x_0,x_1).\end{eqnarray*} Since
\begin{displaymath}
\lim_{n \to \infty}k_i^n = 0,
\end{displaymath}it follows that there exists $N(i,\varepsilon)$
such that$d_i(x_n,x_m) < \varepsilon$ and hence $(x_n,x_m) \in
\mathcal U$ for all $n,m \ge N(i,\varepsilon).$ therefore the
sequence $\{x_n\}$ is a Cauchy sequence in the $d_i$-uniformity on
$X.$\newline let $S_p= \{x_n : n \ge p\}$ for all positive integer
$p$ and let $\beta$ be the filter basis $\{S_p : p=1,2,...\}.$ Then
since $\{x_n\}$ is a $d_i$-cauchy sequence for each $i \in I,$ it is
easyto see that the filter basis $\beta$ is a Cauchy filter in the
uniform space $(X,U).$ To see this we first note that the family
$\beta =\{V^*(i,\varepsilon):i\in I,\}$ is a base for $U$ as $p^*
=\{d_i : i \in I\}.$ Now since $\{x_n\}$ is a $d_i$-Cauchy sequence
in $X,$ there exists a positive integer $p$ such that
 $d_i(x_n,x_m)\le \varepsilon$ for $m\ge p,n\ge p.$ This implies that
$S_p\times S_p\subseteq V^*(i,\varepsilon).$ Thus given any
$\mathcal U\in U,$ we can find an $S_p\in\beta$ such that$S_p\times
S_p\subset \mathcal U.$ Hence $\beta$ is a Cauchy filter in $(X,U).$
Since $(X,U)$ is a $F$-orbitally complete and Hausdorff space,
$S_p\rightarrow z$ for some $z\in X.$ consequently
$F(S_p)\rightarrow Fz$ (follows from the continuity of $F$). Also
$$S_{p+1}\subseteq F(S_p)=\cup \{Fx_n : n\ge P\}$$ for $p=1,2,....$
it follows that $z\in Fz.$ Hence $z$ is a fixed point of $F.$ This
completes our task.$\square$\newline If we take $r = 1$ then the
previous theorem becomes as the followwing theorem.
\begin{theorem} Let $(x,U)$ be a \emph{F-orbitally complete} Hausdorff uniform
space defined by $\{d_i:i\in I\}=p^*.$ and $(2^X,U^*)$ a hyperspace
and let $F:X \rightarrow 2^X$ be a continuous function with $Fx$
compact for each $x \in X$. Assume that
$$min\{H_i(Fx,Fy),d_i(x,Fx)d_i(y,Fy),d_i(y,Fy)\}$$ $$+ a_{i} min\{d_i(x,Fy),d_i(y,Fx)\}\le b_i d_i(x,Fx)+c_i d_i(x,y)$$
for all $i \in I$ and $x,y \in X,$ where $a_i,$ $b_i,$ $c_i$ are
real numbers such that $0<b_i+c_i<1$, then $F$ has a fixed
point.\end{theorem} We denote that if $F$ is a single valued
function on $X$, then we can write $d_i(Fx,Fy) = H_i(Fx,Fy),x,y \in
X, i\in I.$ \newline Thus we obtain the following theorem as a
consequence of the second theorem where we assume $r = 1.$
\begin{theorem} Let $(x,U)$ be a \emph{T-orbitally complete} Hausdorff uniform
space and let $T:X\rightarrow X$ be a \emph{T-orbitally continuous}
function satisfying
$$min\{d_i(Tx,Ty),d_i(x,Tx)d_i(y,Ty),d_i(y,Fy)\}$$ $$+ a_{i} min\{d_i(x,Ty),d_i(y,Tx)\}\le b_i d_i(x,Tx)+c_i d_i(x,y)$$
for all $i \in I$ and $x,y \in X,$ where $a_i,$ $b_i,$ $c_i$ are
real numbers such that $0<b_i+c_i<1$, then $T$ has a \emph{unique
fixed point} whenever $a_i> c_i > 0.$\end{theorem}
\emph{proof}
Define a function $F$ of $X$ into $2^X$ by putting $FX = \{Tx\}$ for
all $x\in X.$ it follows that $F$ satisfies the conditions of second
theorem. Hence $T$ has a fixed point.\newline Now if $a_i > c_i >
0,$ we show that $T$ has a unique fixed point . Assume that $T$ has
two fixed points $Z$ and $w$ which are distinct. Since $d_i(z,Tz) =
0$ and $d_i(w,Tw) = 0,$ then by the condition of the previous
theorem,
$$a_i\min\{d_i(z,Tw),d_i(w,Tz)\}\le c_id_i(z,w)$$ or $$a_id_i(z,w)
\le c_id_i(z,w),$$ $$d_i(z,w) \le \frac {c_i}{a_i} \,d_i(z,w)$$
which is impossible . Thus if $a_i > c_i > 0,$ then $T$ has a unique
fixed point in $X.$ And that is what we wanted.$\square$
\newline \large {COROLLARY}
\newline  \emph{Let} $T$ \emph{be an orbitally continuous self-map of a
$T$-orbitally complete uniform space $(X,U)$ satisfying the
condition
$$\min\{d_i(Tx,Ty),d_i(x,Tx),d_i(y,Ty)\}$$ $$-\min\{d_i(x,Ty),d_i(y,Tx)\}\le
b_id_i(x,Tx) + c_id_i(x,y),$$ $x , y \in X, i\in I$ and
$0<b_i+c_i<1.$ Then for each $x\in X,$ the sequence $\{T^{n}x\}$
converges to a fixed point of $T.$}
\newline
I want to remark that this theorem came to my mind as I thought that
the result obtained by Dhage[2]can be extended with more
generalisation from  metric space to spaces of arbitrary uniform
structure.It is obvious that Theorem 1 and Corollary 1 of his paper
will follow from our last theorem and the corolllary if we replace
the uniform space $(X,U)$ by a metric space.

\end{document}